\documentclass{article}

\usepackage{amsthm,amssymb,latexsym,amsmath, fancyhdr}

\newcommand{\version}{Version 3.0,\ \  January 14, 2011}
\newcommand{\R}{\mathbb{R}} \newcommand{\C}{\mathbb{C}}
 \newcommand{\hh}{\mathbb{H}}
\newcommand{\p}[1]{{\mathbb{P}^{#1}}} 

\newcommand{\pn}{{\mathbb{P}^n}} 
\newcommand{\pd}{{\mathbb{P}^d}} \newcommand{\opd}{{\cal O}_{\mathbb{P}^d}}
 \newcommand{\ol}{{\cal O}_{\ell}}
 \newcommand{\evp}{{\rm ev}_p}
\newcommand{\im}{{\rm Im}~} 
\newcommand{\h}[1]{\otimes H^0(\mathbb{P}^d,{\cal O}_{\mathbb{P}^d}({#1}))}

\newcommand{\calm}{{\cal M}} \newcommand{\cals}{{\cal S}}

\newtheorem{theorem}{Theorem}
\newtheorem{proposition}[theorem]{Proposition}
\newtheorem{lemma}[theorem]{Lemma}
\newtheorem{corollary}[theorem]{Corollary}
\newtheorem{remark}[theorem]{Remark}

\newtheorem{definition}[theorem]{Definition}
\newtheorem{conjecture}[theorem]{Conjecture}
\newtheorem{claim}[theorem]{Claim}

\makeatletter

\@addtoreset{equation}{section}
\@addtoreset{theorem}{section}
\makeatother

\def\Definition{\begin{definition}}
\def\ed{\end{definition}}
\def\Remark{\begin{remark}}
\def\er{\end{remark}}
\def\Claim{\begin{claim}}
\def\ec{\end{claim}}
\def\Lemma{\begin{lemma}}
\def\el{\end{lemma}}
\def\Theorem{\begin{theorem}}
\def\et{\end{theorem}}
\def\Proposition{\begin{proposition}}
\def\ep{\end{proposition}}

\def\goth{\mathfrak}
\def\6{\partial}
\def\endproof{\hbox{\vrule width 4pt height 4pt depth 0pt}}
\newcommand{\restrict}[1]{{\left|_{{\phantom{|}\!\!}_{#1}}\right.}}
\newcommand{\arrow}{{\:\longrightarrow\:}}
\newcommand{\Tw}{\operatorname{Tw}}
\newcommand{\Sec}{\operatorname{Sec}}
\newcommand{\Sym}{\operatorname{Sym}}
\newcommand{\End}{\operatorname{End}}

\begin{document}

\title{Moduli spaces of framed instanton bundles \\ on $\C\p3$ 
and twistor sections of moduli spaces of instantons on $\C^2$}
\author{Marcos Jardim \\ IMECC - UNICAMP \\
Departamento de Matem\'atica \\ Caixa Postal 6065 \\
13083-970 Campinas-SP, Brazil  \\[4mm] Misha Verbitsky \\  
Laboratory of Algebraic Geometry, GU-HSE, \\
7 Vavilova Str., Moscow, Russia, 117312
}

\maketitle

\begin{abstract}
We show that the moduli space $\calm$ of framed instanton
bundles on $\C\p3$ is isomorphic (as a complex manifold)
to a subvariety in the moduli of rational curves of the
twistor space of the moduli space of framed instantons on
$\R^4$. We then use this characterization to prove that
$\calm$ is equipped with a torsion-free affine connection
with holonomy in $Sp(2n,\C)$.
\end{abstract}

\tableofcontents



\section{Introduction}



The operation of complexifying a real algebraic
variety is well known. It is a special case of
a procedure known as ``extension of scalars'',
and is given essentially by tensoring of all
the relevant rings with $\C$. In the converse 
direction, there is an operation called
``Weil's restriction of scalars'', producing
a real algebraic variety (with the same 
topological space) from a complex one.
This operation is given by taking 
real and imaginary parts of all 
algebraic functions.

A composition of these two operations
is called {\bf complexification}: given
a complex algebraic variety of dimension $k$,
it produces another complex algebraic variety, of dimension $2k$.

In complex analytic category, this operation 
is defined only locally. Starting with a 
complex analytic $n$-manifold $X$, it produces
another complex analytic manifold $X_\C$ of dimension $2n$, equipped with an
anticomplex involution $\iota$, such that the 
fixed point set $X_\R:=St(\iota)$ of $\iota$ is equivalent to
$X$ as a real analytic manifold. The manifold $X_\C$ is
defined non-uniquely; only the germ of a neighborhood
of $X_\R$ in $X_\C$ is unique. One of the
ways to produce a complexification is to take
$X_\C:= X \times \bar X$, where $\bar X$ is the
``complex conjugate'' of $\C$, that is,
the same manifold with complex conjugate action of $\C$.

Recall that {\bf a hyperk\"ahler manifold} is
a Riemannian manifold $M$ with a triple
of complex structures $I, J, K\in \End(TM)$
which are K\"aler and satisfy quaternionic
relations. As shown by D. Kaledin \cite{_Kaledin:book_}
and B. Feix \cite{_Feix_},
a complexification of a K\"ahler manifold
is equipped with a canonical hyperk\"ahler structure in a
neighborhood of the fixed point set.
To be more precise, Kaledin and Feix (independently)
constructed a hyperk\"ahler structure on
a neighborhood of the zero section in
the cotangent space $T^*M$. It is also easy to show
that this cotangent space is naturally isomorphic to the
complexification of $M$. This claim is implicit
in Kaledin's work, and found in explicit
form in Feix's Ph. D. thesis, 
\cite{_Feix:thesis_}, Section 2.2.2.

Therefore, it is natural to ask what 
happens if one complexifies a hyperk\"ahler manifold.
Some results in this direction are known.

Recall that {\bf a twistor space} $\Tw(M)$ of a hyperk\"ahler
manifold $M$ is  $\C\p1 \times M$ equipped
with a complex structure which is defined
as follows. Embed the sphere 
$S^2 = \C\p1 \subset \hh$ 
into the quaternion algebra $\hh$ 
as the set of all quaternions $J$ 
with $J^2 = -1$. For every point
$x = \{m\} \times \{J\} \in M \times S^2$, the tangent space $T_x\Tw(M)$ is
canonically decomposed as $T_xX = T_mM \oplus T_J\C\p1$. 
Let $I_J:T_J\C P1 \to T_J\C P1$ be the usual
complex structure operator, and $I_m:T_mM \to T_mM$ 
be the complex structure on $M$ induced by $J \in S^2\subset \hh$.

The operator $I_x = I_m \oplus I_J:T_x\Tw(M) \to T_x\Tw(M)$ satisfies $I_x \circ I_x =
-1$. Moreover, it depends smoothly on the point $x$, hence it defines an almost complex
structure on $\Tw(M)$. This almost complex structure is known to be integrable
(see \cite{_Salamon_}). 

It is well-known that a complexification of a hyperk\"ahler manifold is 
naturally identified with a component called $\Sec(M)$ 
in the moduli of rational curves on its twistor space
(see e.g. \cite{_Verbitsky:hypercomple_} and Subsection \ref{_ratcurves_Subsection_} below).
In fact, it is possible to characterize
hyperk\"ahler manifolds in terms of the
geometric structures on the space $\Sec(M)$.
This approach is quite useful, because
(as suggested by Deligne and Simpson)
it allows to define singular hyperk\"ahler 
varieties (see \cite{_Verbitsky:hypercomple_}
for a precise definition and a 
desingularization theorem).

In this paper we study the geometry 
of $\Sec(M)$ in some detail, and apply
this to obtain a description of geometry
of the space of mathematical instantons
on $\C\p3$.

More precisely, we show that there exists a 
family of integrable foliations $S_v \subset T\Sec(M)$,
parametrized by $v\in \C\p1$, where 
$\dim S_v = \frac 1 2 \dim \Sec(M)$.
For each $C \in \Sec(M)$, the tangent space
$T_C \Sec(M)$ is naturally identified with
the space of sections of $NC$, where
$C$ is identified with $\C\p1$, and
its normal bundle $NC$ (in generic $C$)
with a direct sum of several copies
of ${\cal O}(1)$. We define
$S_v \restrict C$ as the space of all
$x \in NC$ vanishing at $v\in \C\p1 =C$.

This family of foliations form what is known
as a \emph{3-web}, see Subsection \ref{_3_webs_Subsection_} below for a precise definition.
Theory of 3-webs was developed in 1930-ies by
S. S. Chern (\cite{_Chern:3-webs_}), who chose it as a subject of his
doctoral dissertation under Blaschke.
Among other things, Chern constructed a
natural torsion-free connection on a manifold with
a 3-web, called {\bf the Chern connection}.

The 3-webs occurring on $\Sec(M)$ are of 
a particular kind, which we call an $SL(2)$-\emph{web}.
For such an $SL(2)$-web, we obtain 
that the Chern connection is flat on the leaves
of the foliation $S_v\subset T\Sec(M)$, and
its holonomy lies in the centralizer of
$\goth{sl}(2, \C) = \goth{su}(2)\otimes \C$
acting on $\C^{4n}= ({\C}^2)^{2n}$,
as a direct sum of $2n$ weight 1 
representations. As a consequence of Schur's 
lemma, this centralizer is isomorphic $GL(2n, \C)$.

On the manifold $\Sec(M)$, which is
a complexification of a hyperk\"ahler
manifold $M$, it is easy to see
that the Chern connection is in fact a 
complexification of the Levi-Civita
connection $\nabla_{LC}$ on $M$ (it follows immediately
from the uniqueness of the Chern 
connection on 3-webs). Since the 
holonomy of $\nabla_{LC}$ is $Sp(n)$, 
it follows that the holonomy of the Chern connection
lies in its complexification $Sp(n, \C)$.

Recall now that the de Rham algebra of a hyperk\"ahler manifold $M$
is equipped with a natural multiplicative action
of the group $SU(2)$ of unitary quaternions.
A connection (not necessarily Hermitian)
on a vector bundle $B$ is called 
{\bf autodual}, or {\bf NHYM-autodual} if its curvature is
$SU(2)$-invariant as an $\End(B)$-valued 2-form.

Now let $\sigma:\; \Tw(M)\arrow M$ denote the natural projection. It is well known that 
a pullback to the twistor space of an autodual connection on $B\to M$ defines a holomorphic
structure on $\sigma^* B$ (see \ref{_pullba_holo_Lemma_} below and \cite[Lemma 5.1]{_NHYM_}).
In \cite{_NHYM_}, the converse statement was proved (see \ref{_twi_NHYM_Theorem_} below):
a holomorphic bundle on $\Tw(M)$ which is trivial on the rational curves of
form $\{m\} \times\C\p1 \subset M \times\C\p1=\Tw(M)$ is obtained as a pullback of a NHYM-autodual
bundle $(B, \nabla)$ on $M$, defined uniquely.

When $M$ is compact, more can be said. 
Each of the holomorphic bundles on $\Tw(M)$
obtained this way corresponds uniquely
to a section of the twistor projection
$\Tw(W) \arrow \C\p1$, where $W$
is a connected component of the moduli space of holomorphic
bundles on $M$, equipped with a natural
hyperk\"ahler structure (\cite{_Verbitsky:Hyperholo_bundles_}).

One of the main goals of this paper is to show that a
similar correspondence can be defined for the case
$M=\R^4$, and $W$ a connected component of the moduli space of framed 
instantons. More precisely, we prove that the space of holomorphic sections
of the twistor fibration $\Tw(W)\arrow\C\p1$ is biholomorphically equivalent
to moduli space of instanton bundles on $\C\p3$ framed at a line, thus we
conclude that this moduli space has the structure of a $SL(2)$-\emph{web}. 

This result can be used to obtain many kinds of
geometrical information about the instanton spaces.
As a quick application, we show that the moduli
space of instanton bundles on $\C P^3$ framed at a 
given complex line has no complex subvarieties.
Also, we prove that at least one irreducible component 
of this moduli space has the expected dimension
(this result is already known from
\cite{_Tikhomirov:obscure_}).

\hfill

\noindent{\bf Acknowledgments.} 
The first named author is partially supported by the CNPq
grant number 305464/2007-8 and the FAPESP grant number
2005/04558-0. The second named author was partially supported by the 
FAPESP grant number 2009/12576-9, RFBR grant 09-01-00242-a,
 RFBR grant 10-01-93113-NCNIL-a, Science Foundation of 
the SU-HSE award No. 10-09-0015, and AG Laboratory 
GU-HSE, RF government grant, ag. 11 11.G34.31.0023.


\section{Holomorphic 3-webs on complex manifolds}


\subsection{$SL(2)$-webs and Chern connection}
\label{_3_webs_Subsection_}

The following notion is based on
a classical notion of a 3-web, developed
in the 1930-ies by Blaschke and Chern, and much
studied since then.

\Definition
Let $M$ be a complex manifold,
$\dim_\C M = 2n$, and $S_t \subset TM$
a family of $n$-dimensional holomorphic sub-bundles,
parametrized by $t\in \C\p1$.
This family is called {\bf a holomorphic $SL(2)$-web}
if the following conditions are satisfied
\begin{description}
\item[(i)] Each $S_t$ is involutive (integrable), that is, 
$[S_t, S_t] \subset  S_t$.
\item[(ii)] For any distinct points $t, t'\in \C\p1$,
the foliations $S_t$, $S_{t'}$ are transversal:
$S_t \cap S_{t'}=\emptyset$.
\item[(iii)] Let $P_{t,t'}:\; TM \arrow S_t$ be a
  projection of $TM$ to $S_t$ along $S_{t'}$.
Then $P_{t,t'}\in \End(TM)$ generate a 4-dimensional
sub-bundle in $\End(TM)$.
\end{description}
\ed

\Remark
The classical definition of 3-webs 
(see e.g. \cite{_Nagy:3-webs_}) is quite similar:
one is  given three integrable foliations $S_0$, $S_1$ and
$S_\infty$ which are pairwise transversal. 
The $SL(2)$-webs defined above can be obtained
as a special case of a 3-web.
\er

\Remark
The operators $P_{t,t'}\subset \End(M)$
generate a Lie algebra isomorphic to $\mathfrak{gl}(2)$
\er

\Definition
(\cite{_Atiyah:conne_})
Let $B$ be a holomorphic vector bundle
over a complex manifold $M$. A {\bf holomorphic connection}
on $B$ is a holomorphic differential operator
$\nabla:\; B \arrow B \otimes \Omega^1 M$
satisfying $\nabla(fb) = b \otimes df + f \nabla(b)$,
for any holomorphic function $f$ on $M$.
\ed

\Remark
Let $\nabla$ be a holomorphic connection on a holomorphic
bundle, considered as a map $\nabla:\; B \arrow B \otimes \Lambda^{1,0} M$,
and $\bar \6:\; B \arrow B \otimes \Lambda^{0,1} M$ the
holomorphic structure operator. The sum 
$\nabla_f:=\nabla+ \bar\6$ is clearly a connection.
Since $\nabla$ is holomorphic, $\nabla\bar\6 + \bar\6 \nabla=0$,
hence the curvature $\nabla_f^2$ is of type $(2,0)$.
The converse is also true: a $(1,0)$-part of a 
connection with curvature of type $(2,0)$
is always a holomorphic connection.
\er

The following claim is well known (in the smooth setting).
Its holomorphic version is no different.

\Claim\label{_Chern_conne_Claim_}
(See \cite{_Nagy:3-webs_}, Theorem 3.2)
Let $S_1, S_2, S_3$ be a holomorphic 3-web on a
complex manifold $M$. Then there exists a unique
holomorphic connection $\nabla$ on $M$ which 
preserves the foliations $S_i$, and such that its
torsion $T$ satisfies $T(S_1, S_2)=0$.
\ec

\Definition
This connection is called {\bf the Chern connection}
of a 3-web.  It was constructed by Chern in 1936
in his doctoral dissertation under direction of Blaschke 
(\cite{_Chern:3-webs_}).
\ed

\Remark
Please notice that this definition is not
symmetric on $S_1, S_2, S_3$.
\er

\Proposition\label{_torsion-free_Proposition_}
Let $S_t, t\in \C\p1$ be an $SL(2)$-web, and 
$\nabla$ the Chern connection associated with
$S_0, S_1, S_\infty$. Then $\nabla$ is torsion-free.
\ep

{\bf Proof:} 
Let $S_t, t\in \C\p1$ be an $SL(2)$-web, and 
$\nabla$ its Chern connection associated with
$S_0, S_1, S_\infty$.
Clearly, the torsion $T:\; TM \times TM\arrow TM$
preserves the sub-bundles $S_{t}\subset TM$. Indeed, these
bundles are integrable, and $\nabla$ preserves the projection
maps $P_{t,t'}$. Consider a decomposition 
$\Lambda^2 TM= \Lambda^2_{inv} TM\oplus \Lambda^2_+ TM$
onto its invariant part and weight 2 part under the
 $\mathfrak{sl}(2)$-action. From a computation relating the
action of $\mathfrak{sl}(2)$ on $TM$ and the condition
$T(S_t, S_t)\subset S_t$ we infer that $T$ vanishes
on $\Lambda^2_+ TM$ (by Lemma \ref{_Lambda_+_T_Lemma_} below). 
However, an $\mathfrak{sl}(2)$-invariant
form which vanishes on $S_0 \otimes S_1$ must be zero,
by Lemma \ref{_proje_sl_2_2-form_Lemma_} below.
\endproof

\Lemma \label{_Lambda_+_T_Lemma_}
Let $V$ be a weight 1 representation of $\goth{sl}(2)$,
$R$ the set of all $r\in \goth{sl}(2)$
which satisfy $\dim (r V) = \frac 1 2 \dim
V$\footnote{This is equivalent to $q(r,r)=0$, where
$q$ is a Killing form on $\goth{sl}(2)$.}, and
$T:\; \Lambda^2 V \arrow V$ a linear map which 
satisfies
\begin{equation}\label{_respect_proje_Equation}
T(\im r,\im r)\subset \im r, \ \  \forall r\in R.
\end{equation}
Then $T\restrict{\Lambda^2_+V}=0$, where 
$\Lambda^2 V = \Lambda^2_{inv} V \oplus \Lambda^2_+(V)$
is the decomposition of $\Lambda^2 V$ onto the
$\goth{sl}(2)$-invariant and weight 2 component. 
\el

{\bf Proof:} 
Let $x,y$ be a standard basis in an irreducible 
weight 1 representation $H$ of $\goth{sl}(2)=\langle f,g,h\rangle$,
with $h(x)=x, h(y)=-y, f(x)=y, g(y)=x$, and $\Sym^2 H$
be its symmetric square. As a representation of
$\goth{sl}(2)$, the space 
$\Lambda^2_+(V)$ is isomorphic to a direct sum
of several copies of $\Sym^2 H$. On each of these
summands, $T$ induces a map $T_0:\; \Sym^2 H\arrow H$
satisfying \eqref{_respect_proje_Equation}.
Let $T_0^*:\; H \arrow \Sym^2 H$ be its dual
map (with respect to a natural $\goth{sl}(2)$-invariant
pairing). The condition \eqref{_respect_proje_Equation}
now implies that $T_0^*(ax +by)$ is proportional
to $(ax + by)^2\in \Sym^2 H$, for any $a, b\in \C$.
The set of vectors of form $\lambda (ax + by)^2$
is a quadric in $\Sym^2 H$. Since the quadric
contains no 2-dimensional planes, the image
of $T^*_0$ is one-dimensional unless
$T^*_0=0$. This is impossible, because
$T_0^*(ax +by)$ is proportional to $(ax + by)^2$
for all $a, b\in \C$. We proved that $T^*_0=0$.
\endproof

\Lemma \label{_proje_sl_2_2-form_Lemma_}
Let $V$ be a representation of $\goth{sl}(2)$ of weight 1,
$R$ the set of all $r\in \goth{sl}(2)$
which satisfy $\dim (r V) = \frac 1 2 \dim
V$,
$x, y \in R$, and $\xi \in \Lambda^2 V$ an 
$\goth{sl}(2)$-invariant 2-form which
vanishes on $S_x \otimes S_y$, where
$S_x = xV$, $S_y= yV$. Then $\xi=0$.
\el

{\bf Proof:}
It is easy to see that the projectivization 
$\mathbb{P}R$ is isomorphic to $\C\p1$, and the
adjoint action of $SL(2)$ on $\goth{sl}(2)$
induces the standard action of $SL(2)$ on 
$\mathbb{P} R=\C\p1$. Since this action
is bitransitive, the pair $S_x, S_y$
can be transformed to $S_x, S_z$
for any $z\in R$ distinct from $x$,
by an appropriate action of $SL(2)$.
Since $\xi$ is $SL(2)$ invariant
and vanishes on $S_x \otimes S_y$,
it also vanishes on $S_x \otimes S_z$,
for all $z\in R$. Since a sum of all
$S_z$ is $V$, this implies that
$\xi$ vanishes on $S_x \otimes V$.
\endproof

\Remark 
From Proposition \ref{_torsion-free_Proposition_}
it follows immediately that the Chern connection
on a manifold with $SL(2)$-web does not depend
from the choice of three points
$S_0, S_1, S_\infty$ in $\C\p1$. Indeed,
the torsion-free connection
preserving $P_{t,t'}\subset \End(M)$
is unique by Claim \ref{_Chern_conne_Claim_}.
\er

\Theorem
Let $M$ be a manifold equipped with a holomorphic 
$SL(2)$-web. Then its Chern connection is a torsion-free
affine holomorphic connection with holonomy in
$GL(n, \C)$ acting on $\C^{2n}$ as a centralizer of 
an $SL(2)$-action, where $\C^{2n}$ is a direct
sum of $n$ irreducible $SL(2)$-representations
of weight 1. Conversely, every connection with
such holonomy preserves a holomorphic $SL(2)$-web. 
\et

{\bf Proof:}
Let $\nabla$ be a Chern connection associated
with a holomorphic $SL(2)$-web. Then $\nabla$
is holomorphic and commutes with an action
of $\mathfrak{gl}(2)$ generated by the projection
maps $P_{t,t'}\subset \End(M)$. Conversely,
for every such connection, 
$Hol(\nabla) \subset GL(n, \C)$
implies that $\nabla$ preserves a family
of $n$-dimensional sub-bundles of $TM$
parametrized by $\C\p1$. These sub-bundles
are integrable because $\nabla$ is torsion-free.
\endproof

\subsection{An example: rational curves on a twistor
  space}
\label{_ratcurves_Subsection_}

The basic example of holomorphic $SL(2)$-webs
comes from hyperk\"ahler geometry. Let $M$ be a hyperk\"ahler manifold, and 
$\Tw(M)$ its twistor space. Denote by $\Sec(M)$ the
space of holomorphic sections of the
twistor fibration ${\rm Tw}(M)\stackrel\pi\arrow\C\p1$.

We consider $\Sec(M)$ as a complex variety, with the
complex structure induced from the Douady space of rational
curves on $\Tw(M)$. Clearly, for any $C \in \Sec(M)$,
$T_C \Sec(M)$ is a subspace in the space of sections of 
the normal bundle $NC$. This normal bundle is
naturally identified with $T_\pi\Tw(M)\restrict C$, where
$T_\pi\Tw(M)$ denotes the vertical tangent space.

For each point $m \in M$, one has a horizontal section
$C_m:=\{m\} \times \C\p1$ of $\pi$. The space of
horizontal sections of $\pi$ is denoted $\Sec_{hor}(M)$;
it is naturally identified with $M$. It is easy to check
that $N C_m= {\cal O}(1)^{\dim M}$,
hence some neighborhood of $\Sec_{hor}(M)\subset 
\Sec(M)$ is a smooth manifold of
dimension $2\dim M$. It is easy to see
that $\Sec(M)$ is a complexification of
$M = \Sec_{hor}(M)$ considered as a
real analytic manifold (see \cite{_Verbitsky:hypercomple_}). 

Let $\Sec_0(M)$ be a part consisting
of all rational curves $C\in \Sec(M)$
with $NC = {\cal O}(1)^{\dim M}$.
Clearly, $\Sec_0(M)$ is a smooth, open
subvariety in $\Sec(M)$.

On $\Sec_0(M)$, there is an $SL(2)$-web
constructed as follows. For each 
$C \in \Sec_0(M)$ and $t\in \C\p1=C$, 
define $S_t\subset T_C\Sec_0(M)= \Gamma_C (NC)$
as a space of all sections of $NC$ vanishing at $t\in C$.

It is not difficult to check that 
this is an $SL(2)$-web. Transversality
of $S_t$ and $S_{t'}$ is obvious, because a section
of ${\cal O}(1)$ vanishing at two points is zero.
Integrability of $S_t$ is also obvious, because
the leaves of $S_t$ are fibers of the evaluation
map $ev_t:\; \Sec(M) \arrow (M, t)$, mapping $C: \; \C\p1
\arrow \Tw(M)$ to $C(t)$. The last condition is
clear, because $\Gamma_{\C\p1}(V \otimes_\C {\cal O}(1))
= V \otimes_\C \C^2$, and the projection maps $P_{t, t'}$
act on $V \otimes_\C \C^2$ through the second component.



\section{Autodual NHYM-bundles on quaternionic \\ projective spaces}


\subsection{NHYM autodual bundles on hyperk\"ahler manifolds}

Let $M$ be a hyperk\"ahler manifold; recall that the
quaternionic action on $TM$ naturally induces a 
multiplicative action of $SU(2)$ on $\Lambda^*(M)$.

\Definition\label{_autodual_Definition_}
Let $(B, \nabla)$ be a complex vector bundle 
with connection, not necessarily Hermitian,
on a hyperk\"ahler manifold $M$. 
The connection $\nabla$ is called {\bf 
non-Hermitian Yang-Mills  autodual} (NHYM autodual) if its curvature 
$R \in \Lambda^2(M,\End B)$ is $SU(2)$-invariant.
\ed

\hfill

\Remark
The relation between this notion and the notion
of {\em hyperholomorphic bundles} defined in \cite{_Verbitsky:Hyperholo_bundles_} 
is somewhat intricate. A hyperholomorphic bundle on a compact
hyperk\"ahler manifold is a stable bundle with first two Chern
classes $SU(2)$-invariant. It is known that the
Yang-Mills connection of such a bundle (which
exists by Donaldson-Uhlenbeck-Yau, and is unique) 
has $SU(2)$-invariant curvature,
\cite{_Verbitsky:Hyperholo_bundles_}, hence it is ``NHYM
autodual'' in the sense of the above definition.
Moreover, every stable bundle admitting
a NHYM autodual connection has $SU(2)$-invariant
Chern classes, hence it is hyperholomorphic.
However, the NHYM condition does not necessarily
imply stability. 
\er

\hfill

Autodual connections on $M$ give rise to holomorphic
bundles on the twistor space 
$\Tw(M)$ by means of a construction known as
\emph{twistor transform}. This construction turns out to be
essentially invertible. More precisely, let us first
recall the following result.

\begin{lemma} (\cite[Lemma 5.1]{_NHYM_})
\label{_pullba_holo_Lemma_}
Let $(B, \nabla)$ be a complex vector bundle
with connection on a hyperk\"ahler manifold, and 
$\Tw(M)\stackrel{\sigma}{\rightarrow} M$ the standard projection.
The connection $\nabla$ is NHYM autodual if and only if the connection
$\sigma^*\nabla$ has curvature of Hodge type $(1,1)$. 
\end{lemma}

The holomorphic bundle $(\sigma^*B,
{\sigma^*\nabla}^{0,1})$ is called {\bf the twistor transform} of the autodual bundle $(B,\nabla)$. 

Now let ${\cal B}$ be a holomorphic bundle on $\Tw(M)$. We
say that ${\cal B}$ is {\bf trivial on horizontal curves}
if the restriction of ${\cal B}$
to any $C\in\Sec_{hor}(M)$ is trivial (see Subsection \ref{_ratcurves_Subsection_} for notation).

\Remark
\label{_defo_trivi_hori_Remark_}
Let $M$ be a compact hyperk\"ahler manifold. 
Since a trivial bundle is polystable, a small deformation
of a bundle  which is trivial on $\C\p1$ is again
trivial. Since the set of horizontal curves is identified
with $M$, it follows that it is compact. Therefore,
a small deformation of a bundle which is trivial on
horizontal curves is again trivial on horizontal curves.
\er

\Theorem (\cite[Theorem 5.12]{_NHYM_})
\label{_twi_NHYM_Theorem_}
The twistor transform gives an equivalence between the
category of autodual NHYM-bundles on $M$ and the category
of holomorphic bundles on $\Tw(M)$ which are trivial on
horizontal curves. 
\et

For Hermitian autodual bundles on 
quaternionic-K\"ahler manifolds, an analogue of
Theorem \ref{_twi_NHYM_Theorem_}
is proved by T. Nitta in \cite{_Nitta_}, and the NHYM-version
of his theorem can be obtained in the same way as
\cite[Theorem 5.12]{_NHYM_}. The space $\Sec_{hor}(M)$,
in this case, is the set of all holomorphic lines 
of form $\sigma^{-1}(m)$, where $m\in M$ is a point,
and $\sigma:\; \Tw(M) \arrow M$ a standard
projection, and $\Sec(M)$ the space of all
rational curves in $\Tw(M)$ obtained by
deforming curves in $\Sec_{hor}(M)$.

For $\hh\mathbb{P}^k$, the $k$-dimensional
quaternionic projective space ($k\ge1$), one has
$\Tw(M)=\C\mathbb{P}^{2k+1}$. By Nitta's result mentioned above, there exists a 1-1
correspondence between NHYM autodual
bundles on $\hh\mathbb{P}^k$ and holomorphic bundles on
$\C\mathbb{P}^{2k+1}$ which are trivial on horizontal
curves.

In what follows, we will provide a linear algebraic
description of holomorphic bundles on $\C\mathbb{P}^{2k+1}$
which are framed on a horizontal curve.

\subsection{ADHM description of framed bundles on complex \\ projective spaces}

Let $V$ and $W$ be complex vector spaces, with dimensions $c$ and $r$, respectively. Fix $d\ge 0$, and consider the following data ($k=0,\dots,d$):
$$ A_k,B_{k} \in {\rm End}(V) $$
$$ I_k \in {\rm Hom}(W,V) ~~,~~ J_k \in {\rm Hom}(V,W) ~~. $$
Choose homogeneous coordinates $[z_0:\dots:z_d]$ on $\C\pd$ and define
\begin{equation} \label{st1}
\tilde{A} := A_{0}\otimes z_0 + \cdots + A_{d}\otimes z_d \ \ \ {\rm and} \ \ \
\tilde{B} := B_{0}\otimes z_0 + \cdots + B_{d}\otimes z_d ~~.
\end{equation}
These can be regarded as sections of ${\rm Hom}(V,V)\otimes\opd(1)$. Define also:
\begin{equation} \label{st2}
\tilde{I} = I_0\otimes z_0 + \cdots + I_d\otimes z_d \ \ \ {\rm and} \ \ \
\tilde{J} = J_0\otimes z_0 + \cdots + J_d\otimes z_d ~~.
\end{equation}
Similarly, $\tilde{I}$ can be regarded 
as a section of ${\rm Hom}(W,V)\otimes\opd(1)$, while
$\tilde{J}$ can be regarded as a section of ${\rm
  Hom}(V,W)\otimes\opd(1)$.

A $d$-dimensional ADHM datum is a quadruple $\widetilde{X}=(\widetilde{A},\widetilde{B},\widetilde{I},\widetilde{J})$, and can be thought as a point in the affine space
$\widetilde{\mathbf{B}}_d:=\mathbf{B}\h1$, where 
$$ \mathbf{B} := {\rm End}(V)\oplus{\rm End}(V)\oplus{\rm Hom}(W,V)\oplus{\rm Hom}(V,W) . $$
Given a point $p\in\C\pd$, we have a natural evaluation map $\evp:\widetilde{\mathbf{B}}_d\to\mathbf{B}$; for simplicity, we denote $\widetilde{X}(p):=\evp(\widetilde{X})$; note that $\widetilde{X}(p)$ is a $0$-dimensional ADHM datum.

Recall that a $0$-dimensional ADHM datum is said to be:
\begin{enumerate}
\item {\em stable}, if there is no subspace $S\subsetneq V$ such that
$A(S),B(S)\subset S$ and $I(W)\subset S$;
\item {\em costable}, if there is no nontrivial subspace $S\subset V$ such that
$A(S),B(S)\subset S$ and $S\subset \ker J$;
\item{\em regular}, if it is both stable and costable.
\end{enumerate}

The most relevant definition for this paper is the following.

\begin{definition} \label{m-def}
A datum $\widetilde{X}=(\tilde{A},\tilde{B},\tilde{I},\tilde{J})\in\widetilde{\mathbf{B}}_d$ is said to be {\em globally regular} if $\widetilde{X}_p$ is regular for every $p\in\C\pd$.
\end{definition}

We consider here the following generalization of the ADHM equation, which we call the {\em $d$-dimensional ADHM equation}:
\begin{equation} \label{c4}
[ \tilde{A} , \tilde{B} ] + \tilde{I}\tilde{J} = 0 ~.
\end{equation}

Now consider the following action of the group $G=GL(V)$ on $\widetilde{\mathbf{B}}_d$:
\begin{equation}\label{aa}
g\cdot(A_0,\dots,A_d,B_0,\dots,B_d,I_0,\dots,I_d,J_0,\dots,J_d) =
\end{equation}
$$ (gA_0g^{-1},\dots,gA_dg^{-1},gB_0g^{-1},\dots,gB_dg^{-1},gI_0,\dots,gI_d,
J_0g^{-1},\dots,J_dg^{-1}) ~~. $$

Note that 
$$ g\cdot(\tilde{A},\tilde{B},\tilde{I},\tilde{J}) =
(g\tilde{A}g^{-1},g\tilde{B}g^{-1},g\tilde{I},\tilde{J}g^{-1}) $$
where
$$ g\tilde{A}g^{-1} = gA_{0}g^{-1}\otimes z_0 + \cdots + gA_{d}g^{-1}\otimes z_d $$
and so on. In particular, we have $(g\cdot \widetilde{X})(p)=g\cdot\widetilde{X}(p)$. It is easy to see that such action preserves global regularity and the set of solutions of the $d$-dimensional ADHM equation (\ref{c4}).

\begin{definition}
$\calm_d(r,c)$ denotes the set of globally regular solutions of the $d$-dimensional ADHM equation modulo the action of $G$.
\end{definition}


The corresponding geometric objects are instanton bundles on projective spaces. Recall that an {\em instanton bundle} on $\C\pn$ ($n\geq2$) is a locally free coherent sheaf $E$ on $\C\pn$ with $c_1(E)=0$ satisfying the following cohomological conditions:
\begin{enumerate}
\item[(i)] for $n\geq2$, $H^0(E(-1))=H^n(E(-n))=0$;
\item[(ii)] for $n\geq3$, $H^1(E(-2))=H^{n-1}(E(1-n))=0$;
\item[(iii)] for $n\geq4$, $H^p(E(k))=0$, $2\leq p\leq n-2$ and $\forall k$;
\end{enumerate}
The integer $c=-\chi(E(-1))=h^1(E(-1))=c_2(E)$ is called the {\em charge} of $E$.

Moreover, a locally free coherent sheaf $E$ on $\C\pn$ is said to be of trivial splitting type if there is a line $\ell\subset\C\pn$ such that the restriction $E|_\ell$ is the free sheaf, i.e. $E|_\ell\simeq{\cal O}_{\ell}^{\oplus{\rm rk}E}$. A {\em framing} on $E$ is the choice of an isomorphism $\phi:E|_\ell\to{\cal O}_{\ell}^{\oplus{\rm rk}E}$. A {\em framed bundle} is a pair $(E,\phi)$ consisting of a locally free coherent sheaf $E$ of trivial splitting type and a framing $\phi$. Two framed bundles $(E,\phi)$ and $(E',\phi')$ are isomorphic if there exists a bundle isomorphism $\Psi:E\to E'$ such that $\phi'=\phi\circ(\Psi|_\ell)$.

For $n=2$, every locally free coherent sheaf of trivial splitting type is automatically instanton, since the vanishing condition $(i)$ is satisfied. 

The following result, announced in \cite{J-cr}, establishes the relation between framed instanton bundles and solutions of the $d$-dimensional ADHM equations, generalizing the well-known result due to Donaldson \cite{D1}.

\begin{theorem}
There exists a 1-1 correspondence between equivalence classes of globally regular solutions of the $d$-dimensional ADHM equations and isomorphism classes of instanton bundles on $\C\mathbb{P}^{d+2}$ framed at a fixed line $\ell$, where $\dim W={\rm rk}(E)$ and $\dim V=c_2(E)$.
\end{theorem}

The case $d=0$ is originally due to Donaldson \cite{D1}, who also showed that $\calm_0(r,c)$ is isomorphic to the moduli space of rank $r$ framed instantons (i.e. finite energy connections with anti-self-dual curvature) on $\R^4$ of charge $c$.

The case $d=1$ was considered in detail in \cite{FJ2}. In general, the odd case $d=2k+1$ admits the following differential geometric interpretation.
Fix $\ell\subset\C\mathbb{P}^{2k+1}$ to be a horizontal line. By Remark \ref{_defo_trivi_hori_Remark_}, one may regard the set of (framed) holomorphic bundles on
$\C\mathbb{P}^{2k+1}=\Tw(\hh\mathbb{P}^k)$ which are trivial on horizontal curves as an open subset of $\calm_{2k-1}(r,c)$.

In other words, for each $k\ge1$, there is an open subset
of $\calm_{2k-1}(r,c)$ whose points may be interpreted,
via the twistor correspondence, as autodual NHYM-bundles
on $\hh\mathbb{P}^k$ which are framed at a point.

We will finish this paper by taking a closer look at the case $d=1$.

\subsection{Geometric structures on $\calm_1(r,c)$}

Let $W$ be a component of the moduli space of 
Hermitian autodual bundles on a compact hyperk\"ahler 
manifold $M$. Then $W$ can be identified
with a component of the moduli space of stable holomorphic
bundles on $M$; in this case, $W$ is 
hyperk\"ahler, as follows from
\cite{_Verbitsky:Hyperholo_bundles_}.
Indeed, in  \cite{_Verbitsky:Hyperholo_bundles_}
it was shown that the moduli space of stable holomorphic
structures on a bundle $B$ with $SU(2)$-invariant Chern 
classes $c_1(B)$, $c_2(B)$  is hyperk\"ahler.
The invariance of the Chern classes is clear because
the Chern classes are expressed through the the curvature
using the Gauss-Bonnet formula, and the curvature is
$SU(2)$-invariant.

In this case, one may define its twistor space
$\Tw(W)$. In \cite[Theorem 7.2]{_NHYM_} it was shown that
the space $\Sec(W)$ of twistor sections can be identified
with an open subset of the moduli of holomorphic bundles
on $\Tw(M)$.

In this section we establish an analogue to this result
for framed instantons on euclidean space. As it is also
well-known (see for instance \cite{N}), $\calm_0(r,c)$ has
the structure of a smooth hyperk\"ahler manifold with
$\dim_{\R}\calm_0(r,c)=4rc$; let ${\rm Tw}(\calm_0(r,c))$
denote its twistor space; let also
$\cals(r,c):=\Sec(\calm_0(r,c))$, the set of holomorphic
sections of the twistor fibration
$\Tw(\calm_0(r,c))\to\C\p1$.

\begin{theorem}\label{main}
The moduli space $\calm_1(r,c)$ is biholomorphically
equivalent to the space of sections $\cals(r,c)$.
\end{theorem}

In particular, it follows that $\dim \calm_1(r,c) = 4rc$,
and that the smooth locus of $\calm_1(r,c)$ has the
structure of a holomorphic $SL(2)$-web.

\begin{proof} 
Given a framed instanton bundle on $\C\p3$, let
$\widetilde{X}=(\tilde{A},\tilde{B},\tilde{I},\tilde{J})$
be the associated $1$-dimensional ADHM datum; let
$[\widetilde{X}]$ denote its $G$-orbit. Then define a map
$\sigma:\C\p1\to\calm_0(r,c)$ in the following way:
\[ 
  \sigma(p) = [\widetilde{X}(p)] . 
\]
One easily sees that this map is a well-defined, holomorphic map.

Conversely, giving a (holomorphic) map
$\sigma:\C\p1\to\calm_0(r,c)$ is the same as giving a family
of $0$-dimensional ADHM data $(A(p),B(p),I(p),J(p))$ which
is a regular solution of the ADHM equation at each
$p$. But this is precisely a globally regular solution of
the $1$-dimensional ADHM equation, and therefore a point
of $\calm_1(r,c)$.

This establishes a bijection between the sets
$\calm_1(r,c)$ and $\cals(r,c)$. To check that this 
correspondence is a biholomorphic, we use \cite[Theorem 7.2]{_NHYM_},
where the space of fiberwise stable holomorphic
bundles on a twistor space $\Tw(M)$ for a hyperk\"ahler
manifold $M$ is shown to be biholomorphic with a
space $\Sec(W)$ of twistor sections, associated
with the moduli of stable bundles on $M$. 
However, this result cannot be applied directly, because
in our case $M$ is $\C^2$, which is non-compact,
and $W$ is the moduli of framed instantons. 
Nevertheless,
it is easy to check that the proof of 
\cite[Theorem  7.2]{_NHYM_} carries through to framed instantons
without any changes. In fact, the proof of
\cite[Theorem 7.2]{_NHYM_} is based on an 
equivalence between holomorphic bundles
on a twistor space $\Tw(M)$ and NHYM autodual bundles
on $M$ (\cite[Theorem 5.12]{_NHYM_}), 
which is a local result, valid for non-compact $M$.
\end{proof}

\hfill

As an immediate application of Theorem
\ref{main}, we obtain the following corollary.

\begin{corollary}
Let $\calm_1(r,c)$ be the space of framed  instantons
on $\C\p3$. Then $\calm_1(r,c)$ has no compact complex
subvarieties of positive dimension.
\end{corollary}

{\bf Proof:} In \cite[Proposition 8.15]{_NHYM_},
a strictly plurisubharmonic function was constructed
on $\Sec(M)$ for any hyperk\"ahler manifold $M$.
Then, $\calm_1(r,c)\cong \Sec(\calm_0(r,c))$
admits a strictly plurisubharmonic function.
\endproof

\bigskip

Holomorphic rank $2$ bundles $E\to\C\p3$ with vanishing first Chern class and satisfying
$H^0(E)=H^1(E(-2))=0$ are know in the literature as \emph{mathematical instanton bundles}, see for instance \cite{CTT}. These objects have been intensively studied since the 1980's by various authors, see for instance the survey in the Introduction of \cite{CTT} for several references. In particular, every such bundle is (slope and Gieseker) stable.

Let $\mathcal{I}(c)$ denote the moduli space of such bundles; the following is an important open question (c.f. \cite[Conjecture 1.2]{CTT}):

\begin{conjecture}
$\mathcal{I}(c)$ is an irreducible, non-singular, quasi-projective variety of (complex) dimension $8c-3$.
\end{conjecture}

This conjecture is known to be true for $c\le5$ \cite{CTT}, and each case was originally proved by different sets of authors. 

By \cite[Theorem 3]{J-i}, instanton bundles are precisely those obtained as cohomology of a linear monad. Moreover, \cite[Proposition 11]{J-i} tells us that if $E$ is a rank $n-1$ instanton bundle on $\C\pn$, then $H^0(E)=0$. In particular, it follows that rank $2$
instanton bundles on $\C\p3$ are mathematical instanton bundles.

Moreover, the converse is also true: if $E$ is a mathematical instanton bundle, then there is a (unique up to a scalar) symplectic isomorphism between $E$ and its dual $E^*$; one can then use Serre duality to show that $H^2(E(-2))=H^3(E(-3))=0$, thus $E$ is a rank $2$ instanton bundle.

Therefore, there exists a forgetful map
$\psi:\calm_1(2,c)\to\mathcal{I}_{\ell}(c)$, where
$\mathcal{I}_{\ell}(c)$ is the open subset of
$\mathcal{I}(c)$ consisting of mathematical
instanton bundles restricting trivially to a fixed
$\ell$. The fibers of $\psi$ are the sets of all possible
framings (up to equivalence), thus $\calm_1(2,c)$ becomes
a principal $SL(2)$-bundle over $\mathcal{I}_{\ell}(c)$.

Now it follows from Theorem \ref{main} and the
observations in Section \ref{_ratcurves_Subsection_} that
the irreducible component of $\calm_1(2,c)$ containing
the complexification of $\calm_0(2,c)$
has dimension $2\dim\calm_0(2,c)=8c$; we can
then conclude that at least one component of
$\mathcal{I}_{\ell}(c)$ has dimension $8c-3$. 

The following observation seems to be well-known
 (see e.g. \cite{_Hauzer_Langer_}).

\begin{lemma}
If $E$ is a semistable rank $2$ coherent torsion free
sheaf on $\C\p3$ with $c_1(E)=0$, then there is a line
$\ell$ such that $E|_{\ell}\simeq\ol^{\oplus 2}$.
\end{lemma}

On the other hand, notice that if $E$ is nontrivial, then
there exist a line $\ell'$ such that the restricted sheaf
$E|_{\ell'}$ is nontrivial. So for any given line 
$\ell\subset\C\p3$, there are semistable rank $2$ coherent
torsion free sheaves which are trivial at $\ell$ as well
as sheaves that are nontrivial at $\ell$.

\begin{proof}
The lemma is a direct application of the
famous Grauert-M\"ullich theorem
(\cite{_Grauert_Mullich_}).
The restriction of $E$ to a generic line $\ell\subset\p3$
yields a torsion free (hence locally free) sheaf on
$\ell$. If $E|_{\ell}=\ol(a)\oplus\ol(b)$ with $a\ge b$,
Theorem 3.1 of \cite{M} implies that $0\le a-b\le 1$. But
$a+b=c_1(E)=0$, thus we must have $a=b=0$.
\end{proof}

Now let $\mathcal{G}(c)$ denote the moduli space of
S-equivalence classes of semistable torsion-free sheaves
$E$ of rank $2$ on $\p3$ with $c_1(E)=0$, $c_2(E)=c$ and
$c_3(E)=0$; it is a projective variety. $\mathcal{I}(c)$
can be regarded as the open subset of $\mathcal{G}(c)$
consisting of those locally free sheaves satisfying
$H^1(E(-2))=0$. 

Two important facts follow from our previous lemma. First,
for any fixed line $\ell\subset\p3$, $\mathcal{I}(c)$ is
contained in $\overline{\mathcal{I}_{\ell}(c)}$, where the
closure is taken within $\mathcal{G}(c)$, thus
$\mathcal{I}(c)$ is irreducible if and only if
$\mathcal{I}_{\ell}(c)$ is. Second, $\mathcal{I}(c)$ is
covered by open subsets of the form
$\mathcal{I}_{\ell}(c)$, but it is not contained within
any such sets, thus $\mathcal{I}(c)$ and
$\mathcal{I}_{\ell}(c)$ must have the same dimension, and
one is nonsingular if and only if the other is as well.

Summing up our conclusions, we have proved that the Conjecture
is true if and only if, for some line $\ell\subset\C\p3$,
the quasi-projective variety $\mathcal{I}_\ell(c)$ is
irreducible, non-singular and of dimension
$8c-3$. In particular, we have established that
$\mathcal{I}(c)$ possesses an irreducible component of dimension $8c-3$ for 
all $c\geq 1$; existence of an irreducible component of the expected
dimension was already known, since there are examples several examples of mathematical instanton bundles that are unobstructed, i.e. $H^2({\rm End}(E))=0$; see also \cite{_Tikhomirov:obscure_} for a more elaborate result.

Moreover, $\mathcal{I}_\ell(c)$ is irreducible and
nonsingular if and only if $\calm_1(2,c)$ is. We hope that
the geometric and ADHM-type descriptions of this space
given in this paper will be valuable tools in an attempt
to prove the Conjecture.


\end{document}